\pdfoutput=1
\documentclass[12pt,a4paper]{article}
\usepackage[T1]{fontenc}
\usepackage{amsmath,amssymb,amsthm,mathtools}
\usepackage[left=3.3cm,right=3.3cm,top=3.2cm,bottom=3.5cm]{geometry}
\usepackage{enumitem}
\usepackage{xcolor}
\usepackage{hyperref}
\hypersetup{colorlinks=true,linkcolor=blue!55!black,citecolor=blue!55!black,urlcolor=blue!55!black}
\theoremstyle{plain}
\newtheorem{theorem}{Theorem}[section]
\newtheorem{proposition}[theorem]{Proposition}
\newtheorem{lemma}[theorem]{Lemma}
\newtheorem{corollary}[theorem]{Corollary}
\newcommand{\1}{\mathbf{1}}
\newcommand{\RR}{\mathbb{R}}
\newcommand{\ZZ}{\mathbb{Z}}
\DeclareMathOperator{\sgn}{sgn}
\usepackage{lmodern}
\usepackage{graphicx}
\usepackage{needspace}
\usepackage{float}
\usepackage[labelfont=bf,labelsep=period,margin=1.5em]{caption}
\DeclareCaptionFont{tight}{\small}
\makeatletter
\def\thm@space@setup{\thm@preskip=14pt plus 4pt minus 2pt \thm@postskip=\thm@preskip}
\makeatother
\newcommand{\NN}{\mathbb{N}}
\newcommand{\cP}{\mathcal{P}}
\newcommand{\cG}{\mathcal{G}}
\newcommand{\cC}{\mathcal{C}}
\newcommand{\cD}{\mathcal{D}}
\newcommand{\eps}{\varepsilon}
\title{Discrepancy and numerical\\ integration for the worst measure}
\author{Mario Ullrich\thanks{Department of Quantum Computing \& Institute of Analysis,
Johannes Kepler University Linz, Austria. E-mail: \texttt{mario.ullrich@jku.at}. Date: November 27, 2026.}}
\date{}
\begin{document}\maketitle
\begin{abstract}\noindent
For $d\ge3$ and every large $n$, there is a probability measure on $[0,1)^d$, supported on
$O(n\log n)$ points, such that every set of $n$ points has star
discrepancy at least $(\log n)^{d-1-o(1)}/n$ with respect to it. In fact, we prove a lower
bound for quadrature rules with arbitrary weights, that is, for the worst-case error of numerical
integration with respect to this measure on functions of bounded variation. By a recent upper
bound for the Lebesgue measure, such measures are harder to approximate than the Lebesgue
measure. This answers a question of Aistleitner, Bilyk and Nikolov. We also show bounds
for densities that are bounded in $L_2$ and for the $L_p$-discrepancy. The new idea is to use the small
ball inequality in two forms: to construct many measures, and to separate their distribution
functions.
\end{abstract}

\section{Introduction}

Let $d\ge2$ and $D=[0,1)^d$. Let $\mu$ be a probability measure on $D$, and let
$Q_n(f)=\sum_{i=1}^nw_if(x_i)$ be a quadrature rule with nodes $x_1,\dots,x_n\in D$ and weights
$w_1,\dots,w_n\in\RR$. The star discrepancy
\[
D^*(Q_n;\mu)=\sup_{t\in D}\Bigl|\mu([0,t])-\sum_{i\colon x_i\le t}w_i\Bigr|
\]
measures how well $\mu$ is approximated by the atomic measure with the weights $w_i$ at the
points $x_i$. It is the worst-case error of $Q_n$ for $\int f\,d\mu$ on a unit ball of
functions of bounded variation (Section~\ref{sec:rem}). For a point set $\cP_n=\{x_1,\dots,x_n\}\subset D$
(points may repeat), we write $D^*(\cP_n;\mu)$ for the classical star discrepancy, that is,
for the rule with equal weights $w_i=1/n$.

\medskip\noindent\textbf{Lebesgue measure.} There are point sets with
$D^*(\cP_n;\lambda)\lesssim(\log n)^{d-1}/n$ (Halton \cite{Hal}). Schmidt \cite{Sch} proved the optimal
lower bound $(\log n)/n$ for $d=2$, and Bilyk, Lacey and Vagharshakyan \cite{BLV} the bound
$(\log n)^{\alpha_d}/n$ with some $\alpha_d>\frac{d-1}2$ for $d\ge3$; here $\frac{d-1}2$ is the
exponent for the $L_2$-discrepancy, see below. In a recent preprint,
Wengiel \cite{W} proves the small ball inequality for Haar functions up to a power of a
logarithm, the lower bound $(\log n)^{d/2}/n$ for $d=3$ and up to a power of $\log\log n$ for
$d\ge4$, and constructs point sets with
\begin{equation}\label{eq:Wup}
D^*(\cP_n;\lambda)\ \lesssim\ \frac{(\log n)^{d-1-\gamma_d}}n,\qquad
\gamma_3=\tfrac2{31},\quad\gamma_4=\tfrac1{12},\quad\gamma_d=\tfrac18\ \ (d\ge5),
\end{equation}
see \cite[Theorem~1.6 and Proposition~8.30]{W}.

\medskip\noindent\textbf{General measures.} Beck \cite{Beck} proved that for \emph{every}
$\mu$ there are point sets with $D^*(\cP_n;\mu)\lesssim(\log n)^{2d}/n$. Aistleitner and Dick
\cite{AD14} improved the exponent to $\frac{3d+1}2$, and Aistleitner, Bilyk and Nikolov
\cite{ABN} to $d-\frac12$. In the other direction, Chen
\cite{Chen} extended the lower bound for the $L_p$-discrepancy of the
Lebesgue measure, $1<p<\infty$ (see below), to every measure with a density, also for rules with
arbitrary weights; so no density is easier than the Lebesgue measure, up to constants.
Aistleitner, Bilyk and Nikolov \cite{ABN} asked:
\begin{quote}\itshape
Is there any measure which is more difficult to approximate by finite atomic measures with
equal weights than Lebesgue measure?
\end{quote}
Their intuitive answer was no: it is ``tempting to assume that the
Lebesgue measure is essentially extremal'', since a measure that is difficult to approximate should
be spread out everywhere.
For $d=1$, the Lebesgue measure is the hardest, and some discrete measures are easier, by
Fairchild, Goering and Wei{\ss} \cite{FGW}; see Brandolini, Colzani and Travaglini \cite{BCT} for
sharp results on balls in the torus.

Our main result is the following lower bound for the worst measure. It rests on the small
ball inequality of Wengiel \cite{W}, which also determines $\chi_d$ (Section~\ref{sec:sb}).

\begin{theorem}[Star discrepancy for the worst measure]\label{thm:main}
Let $d\ge3$. There are $c_d>0$, $C_d,\chi_d<\infty$ and $n_d\in\NN$ such that for every
$n\ge n_d$ there is a probability measure $\mu$ on $[0,1)^d$, supported on at most $C_dn\log n$
points, such that every quadrature rule $Q_n$ with $n$ nodes in $[0,1)^d$ and real
weights satisfies
\[
D^*(Q_n;\mu)\ \ge\ c_d\,\frac{(\log n)^{d-1}}{n\,(\log\log n)^{\chi_d}} .
\]
In particular, this holds for every point set $\cP_n\subset[0,1)^d$.
\end{theorem}

\begin{figure}[ht]\centering
\includegraphics[width=0.34\linewidth]{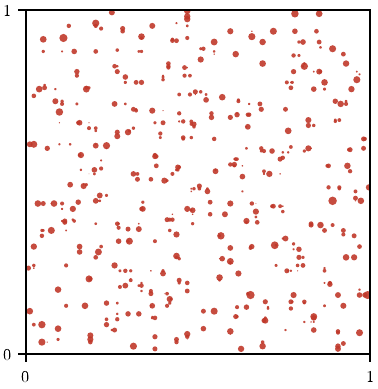}
\caption{A hard measure for $d=2$, supported on about $500$ points (Section~\ref{sec:hard}).
The area of a point is proportional to its weight.}\label{fig:intro}
\end{figure}

\needspace{5\baselineskip}
By the upper bound \eqref{eq:Wup}, for every large~$n$, some
measure is harder to approximate than the Lebesgue measure, even if arbitrary weights are
allowed; the measures depend on~$n$. One measure suffices for infinitely many $n$:

\begin{corollary}[One measure for infinitely many $n$]\label{cor:one}
Let $d\ge3$, $\delta>0$ and let $\gamma_d$ be as in \eqref{eq:Wup}. There is a probability
measure $\mu$ on $[0,1)^d$ with countable support such that
\[
\inf_{Q_n}D^*(Q_n;\mu)\ \ge\ (\log n)^{\gamma_d-\delta}\,\inf_{\cP_n}D^*(\cP_n;\lambda)
\qquad\text{for infinitely many }n,
\]
where the first infimum is over all rules with $n$ nodes and real weights.
\end{corollary}

In this sense, the
answer to the question of Aistleitner, Bilyk and Nikolov is yes; for all large $n$ this is
open. By their upper bound,
the exponent of $\log n$ for the worst measure
lies between $d-1$ and $d-\frac12$, up to powers of $\log\log n$, with or without weights.

\bigskip
For $1\le p<\infty$, the \emph{$L_p$-discrepancy} of $Q_n$ with respect to $\mu$ is
\[
D_p(Q_n;\mu)\ =\ \Bigl(\int_D\Bigl|\mu([0,t])-\sum_{i\colon x_i\le t}w_i\Bigr|^p\,dt\Bigr)^{1/p}.
\]
For $1<p<\infty$, we have the following bounds (Section~\ref{sec:L2}). The upper bound follows from \cite{ABN}, the
lower bound from bounds for entropy numbers of Bilyk, Neumayer and T.~Ullrich \cite{BNU} and
the duality of metric entropy of Artstein, Milman, Szarek and Tomczak-Jaegermann
\cite{AMS,AMST}.

\needspace{10\baselineskip}
\begin{theorem}[$L_p$-discrepancy for the worst measure]\label{thm:L2intro}
Let $d\ge2$ and $1<p<\infty$. There are $c_{d,p}>0$, $C_d<\infty$ and $n_{d,p}\in\NN$ such
that for all $n\ge n_{d,p}$
\[
c_{d,p}\,\frac{(\log n)^{d-\frac32}}n\ \le\ \sup_\mu\,\inf_{Q_n}D_p(Q_n;\mu)
\ \le\ C_d\,\frac{(\log n)^{d-\frac12}}n ,
\]
where the supremum is over all probability measures $\mu$ on $[0,1)^d$ and the infimum over all
quadrature rules $Q_n$ with $n$ nodes in $[0,1)^d$ and real weights. The same bounds hold if
the infimum is taken only over point sets $\cP_n$.
\end{theorem}

For the Lebesgue measure, the order of the $L_p$-discrepancy, $1<p<\infty$, is
$(\log n)^{(d-1)/2}/n$, by Roth \cite{R54,R80}, Schmidt \cite{Sch77} and Chen \cite{Chen80}. So, for $d\ge3$, the worst measure is harder also in $L_p$.
Table~\ref{tab} collects the exponents, also for
densities that are bounded in $L_2$ (Section~\ref{sec:L2}).

\begin{table}[H]\centering\small
\renewcommand{\arraystretch}{1.35}
\begin{tabular}{|l|c|c|c|}\hline
& & \multicolumn{2}{c|}{worst measure}\\
& Lebesgue measure & density with $\|\rho\|_2\le2$ & arbitrary\\\hline
star discrepancy & $\frac d2\le\alpha\le d-1-\gamma_d$ & $d-\frac32\le\alpha\le d-\frac12$ &
$d-1\le\alpha\le d-\frac12$\\
$L_p$-discrepancy & $\alpha=\frac{d-1}2$ & $d-2\le\alpha\le d-1$ &
$d-\frac32\le\alpha\le d-\frac12$\\\hline
\end{tabular}
\caption{The exponent $\alpha$ of $(\log n)^\alpha/n$ in the smallest discrepancy with $n$
points, $d\ge3$, up to powers of $\log\log n$. The
lower bounds hold for real weights, except $\frac d2$ (equal weights). The upper bounds hold for
equal weights, except $d-1$ (weights). In the second row, $1<p<\infty$, except for densities,
where $p=2$.}\label{tab}
\end{table}

\noindent\textbf{Outline.} Section~\ref{sec:idea} sketches the proof of Theorem~\ref{thm:main},
and Sections~\ref{sec:sb}--\ref{sec:proof} contain it. Section~\ref{sec:L2} treats the
$L_p$-discrepancy, and Section~\ref{sec:rem} contains remarks.

\medskip\noindent\textbf{Notation.} For $x,t\in\RR^d$, we write $x\le t$ if
$x_{\ell}\le t_{\ell}$ for all $\ell$, and $x<t$ if $x_{\ell}<t_{\ell}$ for all $\ell$. For
$x,t\in D$, let $[0,t]=\{x\in D\colon x\le t\}$,
$[0,t)=\{x\in D\colon x<t\}$, $[x,1)=\{t\in D\colon x\le t\}$ and $\1=(1,\dots,1)$. Measures are finite signed Borel measures
on $D$, $|\nu|$ is the total
variation measure, $\delta_x$ the point mass at $x$, $\lambda$ the Lebesgue measure, and $\|\cdot\|_p$ the norm of
$L_p=L_p(D,\lambda)$. We write $a_n\lesssim b_n$ or $b_n\gtrsim a_n$ if $a_n\le Cb_n$ for all $n$ with a constant
$C$ independent of $n$, and $a_n\approx b_n$ if $a_n\lesssim b_n\lesssim a_n$. Constants $c,C$ depend only on $d$ unless stated otherwise; their value
may change from line to line.

\needspace{40\baselineskip}
\section{Proof strategy}\label{sec:idea}

We sketch the proof of Theorem~\ref{thm:main}, ignoring powers of $\log\log n$ and constants
that depend only on $d$. For a box $R=I_1\times\dots\times I_d$, let
$\Lambda_R(x)=\prod_\ell\operatorname{dist}(x_\ell,\RR\setminus I_\ell)$ be its tent function;
$\int\Lambda_R\,d\mu$ is the tent mass of $\mu$ in $R$. Let $\cD_m$ be the set of the
$K_m\approx2^mm^{d-1}$ dyadic boxes of volume $2^{-m}$, and let $S(R)\in\cD_m$ be the mirror image
of $R\in\cD_m$ under $x\mapsto\1-x$.

\medskip\noindent\emph{Step 1: many measures} (Section~\ref{sec:hard}). A small ball inequality of
Wengiel \cite{W} for sums of tent functions (Theorem~\ref{thm:tent}) and the separation theorem
for convex sets give a probability measure $\mu_\eta$ for every
$\eta=(\eta_R)_{R\in\cD_m}\in\{-1,1\}^{\cD_m}$ with $\eta_{S(R)}=-\eta_R$, such that
\[
\int\Lambda_R\,d\mu_\eta-\int\Lambda_{S(R)}\,d\mu_\eta\ =\ c_\eta\,\eta_R\quad\text{for all }
R\in\cD_m,\qquad c_\eta\approx m^{-\frac{d-2}2}2^{-2m}
\]
(Lemma~\ref{lem:hard}). For the Lebesgue measure, this difference is zero. There are
$2^{K_m/2}$ such $\eta$, and each $\mu_\eta$ has at most $K_m/2+1$ atoms (Figure~\ref{fig:intro}).

\medskip\noindent\emph{Step 2: their distribution functions are far apart}
(Section~\ref{sec:hard}). By another small ball inequality
(Proposition~\ref{prop:sb}), every signed measure $\nu$ satisfies
\[
\sup_t\,\bigl|\nu([0,t])\bigr|\ \gtrsim\ m^{-\frac{d-2}2}\sum_{R\in\cD_m}\Bigl|\int\Lambda_R\,d\nu\Bigr|
\tag{$*$}
\]
up to an error term. Let $\eta$ and $\eta'$ differ in at least $K_m/4$ entries, and let
$\nu=\mu_\eta-\mu_{\eta'}$. If $\eta_R\ne\eta'_R$, then
$|\int\Lambda_R\,d\nu|+|\int\Lambda_{S(R)}\,d\nu|\ge c_\eta+c_{\eta'}$ by Step~1 (Lemma~\ref{lem:sep}). So the sum in
$(*)$ is $\gtrsim K_mm^{-\frac{d-2}2}2^{-2m}$, and $(*)$ gives
\[
\sup_t\,\bigl|\mu_\eta([0,t])-\mu_{\eta'}([0,t])\bigr|\ \gtrsim\
\underbrace{m^{-\frac{d-2}2}}_{(*)}\cdot\underbrace{K_m}_{\text{pairs}}\cdot
\underbrace{m^{-\frac{d-2}2}2^{-2m}}_{\text{difference}}\ \approx\ m\,2^{-m} .
\]

\medskip\noindent\emph{Step 3: counting} (Section~\ref{sec:proof}). Here is the idea for point
sets, that is, for equal weights $\frac1n$. If the points are rounded to a grid with $n^C$ cells,
there are at most $n^{Cn}$ sets of $n$ points. If a point set $\cP_n$ has star discrepancy much smaller than $m2^{-m}$ with respect to
$\mu_\eta$ and to $\mu_{\eta'}$, then the distribution functions of $\mu_\eta$ and $\mu_{\eta'}$ are
closer than $m2^{-m}$, so $\eta$ and $\eta'$ differ in fewer than $K_m/4$ entries, by Step~2. So
each point set is good only for the $\eta$ in a Hamming ball of radius $K_m/4$, a proportion
$2^{-cK_m}$ of all $\eta$. If $K_m\ge Cn\log n$, then $n^{Cn}2^{-cK_m}$ is small, and for most $\eta$
every point set has $D^*(\cP_n;\mu_\eta)\gtrsim m2^{-m}$ (Lemma~\ref{lem:count}). With
$K_m\approx n\log n$, we get $m2^{-m}\approx(\log n)^{d-1}/n$, the bound of
Theorem~\ref{thm:main}. For real weights, Lemma~\ref{lem:few} replaces the count of point sets:
it approximates the distribution functions of all rules by $2^{Cn\log n}$ functions, whatever
the weights are.

\section{Distribution functions of signed measures}\label{sec:sb}

\noindent\textbf{Dyadic boxes.} A dyadic
interval of depth $j\in\NN_0$ is $I=[k2^{-j},(k+1)2^{-j})\subseteq[0,1)$; $I^-$ and $I^+$ are
its left and its right half, and $h_I=\1_{I^-}-\1_{I^+}$ is the Haar function of $I$. A dyadic box
is $R=I_1\times\dots\times I_d$ with dyadic intervals $I_\ell$. Its shape is
$r(R)=(r_1,\dots,r_d)\in\NN_0^d$, where $r_\ell$ is the depth of $I_\ell$; we write
$|r|=r_1+\dots+r_d$, so that $|R|=2^{-|r(R)|}$. The Haar function of $R$ is
$h_R(x)=\prod_\ell h_{I_\ell}(x_\ell)$. The boxes of one shape partition $D$, and Haar
functions of distinct boxes are orthogonal with $\|h_R\|_2^2=|R|$. Let $\cD_m$ be the set of
dyadic boxes of volume $2^{-m}$ and
\[
K_m=\#\cD_m=2^m\binom{m+d-1}{d-1}\ \ge\ \frac{2^m(m+1)^{d-1}}{(d-1)!}.
\]
Here $\binom{m+d-1}{d-1}$ is the number of shapes $r$ with $|r|=m$, and each shape has
$2^m$ boxes. The \emph{tent} of $I$ is $\Lambda_I(y)=\min\{y-\inf I,\sup I-y\}$ for $y\in I$
and $\Lambda_I(y)=0$ otherwise, and
$\Lambda_R(x)=\prod_{\ell}\Lambda_{I_{\ell}}(x_{\ell})$. Then
$\Lambda_R$ is continuous, vanishes outside $R$, $0\le\Lambda_R\le2^{-d}|R|$, and
$|\Lambda_R(x)-\Lambda_R(y)|\le\sum_{\ell}|x_{\ell}-y_{\ell}|$. We call
$\int\Lambda_R\,d\nu$ the \emph{tent mass} of a measure $\nu$ in~$R$.

The Haar coefficients of a distribution function are tent masses.

\begin{lemma}\label{lem:haar}
Let $\nu$ be a measure.
\begin{enumerate}[label=(\roman*),leftmargin=2em,itemsep=0pt]
\item $\int_D\nu([0,t])\,h_R(t)\,dt=(-1)^d\int\Lambda_R\,d\nu$ for every dyadic box $R$, and
the same holds with $[0,t)$ in place of $[0,t]$.
\item $\sum_{r(R)=s}\bigl|\int\Lambda_R\,d\nu\bigr|\le2^{-d}2^{-|s|}|\nu|(D)$ for every
$s\in\NN_0^d$.
\end{enumerate}
\end{lemma}
\begin{proof}
(i) For $p\in[0,1)$,
$\int_{\{x\ge p\}}h_I(x)\,dx=\int_{\{x>p\}}h_I(x)\,dx=-\Lambda_I(p)$: both sides vanish
for $p\notin I$, and the integral is $\inf I-p$ for $p\in I^-$ and $p-\sup I$ for
$p\in I^+$. Since $\nu([0,t])=\int\1_{\{p\le t\}}\,d\nu(p)$, Fubini's theorem and the product
structure give the claim; for $[0,t)$ we use $\{p<t\}$.
(ii) follows from
$|\int\Lambda_R\,d\nu|\le2^{-d}|R|\,|\nu|(R)$, since the boxes of one shape are disjoint.
\end{proof}

\noindent\textbf{Test measures.} The following theorem is contained in Wengiel's preprint
\cite{W}; see the proof. For given numbers $\sigma_R$, it provides a measure $\omega$ that acts on functions that are constant
on dyadic boxes like a density with the Haar coefficients $\kappa\theta_m\sigma_R$ on $\cD_m$, up
to small terms on smaller boxes. In \eqref{eq:theta} below,
$\zeta_d$ is the exponent of the small ball inequality \cite[Theorem~1.2]{W}.

\begin{theorem}[Test measures]\label{thm:W}
Let $d\ge3$. For $m\in\NN$, let $m^+=m+3d+\lceil2d\log_2(m+2)\rceil$, and let $N_m$ be the
number of shapes $j\in\NN_0^d$ with $|j|\le m^+$. There are $c>0$, $C<\infty$ and
$m_0\in\NN$ with the following property, where
\begin{equation}\label{eq:theta}
\theta_m=c\,(m+1)^{-\frac{d-2}2}(\log(m+2))^{-\zeta_d}.
\end{equation}
Let
$m\ge m_0$, $\kappa\in(0,1]$ and $\sigma_R\in[-1,1]$ for $R\in\cD_m$. Then there is a
measure $\omega\ge0$ on $D$ with $\omega(D)\le1$ such that:
\begin{enumerate}[label=(\alph*),leftmargin=2em,itemsep=0pt]
\item $\omega(R)\le2|R|$ for every dyadic box $R$ with $|r(R)|\le m^+$;
\item for every shape $j$ with $|j|\le m^+$ and every bounded function $g$ that is constant
on the boxes of shape $j$,
\[
\Bigl|\int g\,d\omega-\int g\,\Bigl(1+\sum_{s<j}f_s\Bigr)\,d\lambda\Bigr|\ \le\
\frac{\|g\|_\infty}{2^dN_m},
\]
where $f_s$ is a linear combination of the Haar functions $h_R$ with $r(R)=s$,
\[
f_s=\begin{cases}0&\text{if }|s|<m,\\ \kappa\theta_m\sum_{r(R)=s}\sigma_Rh_R&\text{if }|s|=m,
\end{cases}
\]
and $\|f_s\|_\infty\le C(\kappa\theta_m)^2(|s|-m+1)^{2d}$ if $|s|>m$.
\end{enumerate}
\end{theorem}
\begin{proof}
We refer to the first version of the
preprint \cite{W},
arXiv:2611.24934v1. The
measure $\omega$ is the trimmed measure $\nu'$ of \cite[Lemma~6.2]{W}, built from the binary model
of \cite[Section~4.1]{W} with, in the notation there, shell index $h=m$, lengths $T=m^+$, $C_0=2$, $Q=2d$, the
amplitude $c_*(m+1)^{-\frac{d-2}2}L^{-\zeta_d}$ with $L\approx\log(m+2)$, which is at least
$\theta_m$ if $c$ is small, and the source on the box $R$ equal to $\kappa\sigma_R$ times
$\theta_m$ divided by this amplitude. Part (a) is \cite[(6.3)]{W}. The functions
$1+\sum_{s<j}f_s$ are the densities $W_j$ of \cite[Lemma~6.1]{W}, and the bound for
$\|f_s\|_\infty$ is part of that lemma: the construction depends only on the product of amplitude
and source, here $\kappa\theta_m\sigma_R$ \cite[Section~4.9.2]{W}, and \cite[Lemma~5.1]{W}
holds for every amplitude in $(0,1]$, here $\kappa\theta_m$. Part
(b) is \cite[(6.4)]{W}: there the right-hand side is $\delta\|g\|_\infty$ with
$\delta\le C(m^+)^{d-1}(m+1)^{-2d}$ by \cite[Lemma~6.2, (6.2)]{W}, and
$2^dN_m\delta\le2^d(m^++1)^d\delta\le1$ for large $m$.
\end{proof}

The next proposition is a small ball inequality for the distribution function of a measure
$\nu$: if $\nu$ has tent mass in many boxes of $\cD_m$, then $\nu([0,t])$ is large for some $t$.
For $\nu=\sum_i\delta_{x_i}-n\lambda$, this is the inequality behind Wengiel's lower bound for the
star discrepancy, and the proof follows the one in \cite{W}, with the discrepancy function
replaced by the distribution function of $\nu$.

\begin{proposition}[A small ball inequality for distribution functions]\label{prop:sb}
Let $d\ge3$, and let $\theta_m$ be as in Theorem~\ref{thm:W}. There are $C<\infty$ and
$m_0\in\NN$ such that for every $m\ge m_0$, every $\kappa\in(0,1]$ and every measure~$\nu$,
\[
\sup_{t\in D}|\nu([0,t])|\ \ge\ \frac{\kappa\theta_m}3\sum_{R\in\cD_m}\Bigl|\int\Lambda_R\,d\nu\Bigr|
-C\Bigl(\kappa^2\theta_m^2(m+1)^{d-1}+\frac1{m+1}\Bigr)2^{-m}|\nu|(D).
\]
\end{proposition}

\begin{proof}
Let $V(t)=\nu([0,t))$ for $t\in D$ and $B=\sup_{t\in D}|V(t)|$. Since $[0,t)$ is empty or the
union of an increasing sequence of boxes $[0,s]$, $B$ is at most the left-hand side, and it
suffices to prove the inequality with $B$ on the left. We use half-open boxes because the sets
$\{x\in D\colon x>p\}$ are disjoint unions of dyadic boxes; see the fine part below. By
Lemma~\ref{lem:haar}(i), $\int_DV\,h_R\,d\lambda=(-1)^d\int\Lambda_R\,d\nu$ for every dyadic
box~$R$. We apply Theorem~\ref{thm:W} with $\kappa$ and $\sigma_R=\sgn\int_DV\,h_R\,d\lambda$
(with $\sgn0=1$).

\emph{Projections.} Let $P_k$ be the mean value operator on $[0,1)$ for the dyadic intervals of
depth $k$: the function $P_kf$ is constant on each of these intervals, and equal to the mean value
of $f$ there. Let $\Delta_k=P_k-P_{k-1}$ ($P_{-1}=0$), $\Delta_j=\bigotimes_\ell\Delta_{j_\ell}$ for
$j\in\NN_0^d$, and $\Pi=\sum_{|j|\le m^+}\Delta_j$. The function $\Delta_jf$ is constant on the
boxes of shape $j$, $\|\Delta_jf\|_q\le2^d\|f\|_q$ for $q\in\{1,\infty\}$, and
$\int(\Delta_jf)\,g\,d\lambda=\int f\,\Delta_jg\,d\lambda$ for bounded $f,g$, since the $P_k$ are
orthogonal projections of $L_2[0,1)$. If $f$ is constant on the boxes of a shape $j$, then
$\Delta_{j'}f=0$ unless $j'\le j$, and $\sum_{j'\le j}\Delta_{j'}f=\bigotimes_\ell P_{j_\ell}f=f$;
so $\Pi f=f$ if moreover $|j|\le m^+$. Finally, $\Delta_01=1$ and $\Delta_j1=0$ for $j\ne0$, and for a
dyadic box $R$, $\Delta_jh_R=h_R$ if $j=r(R)+\1$ and $\Delta_jh_R=0$ otherwise.

\emph{The fine part.} We show that
\[
\Bigl|\int(V-\Pi V)\,d\omega\Bigr|\ \le\ \frac{C2^{-m}}{m+1}\,|\nu|(D).
\] For $p\in D$, let $O_p=\{x\in D\colon x>p\}$, so that
$V(x)=\int\1_{O_p}(x)\,d\nu(p)$ for $x\in D$. Since $\Pi$ is an integral operator with a bounded
kernel, Fubini's theorem reduces the claim to
$|\int(\1_{O_p}-\Pi\1_{O_p})\,d\omega|\le C2^{-m}/(m+1)$ for every $p\in D$.

In one variable, let $q\in[0,1)$. A number $y\in[0,1)$ satisfies $y>q$ if and only if $y$ has the
digit one at the first binary digit where $y$ and $q$ differ. For $k\in\NN$, the numbers $y$ for
which this is the $k$-th digit form a dyadic interval $J_k(q)$ of depth $k$, or the empty set,
and $J_k(q)$ lies in the dyadic interval of depth $k-1$ that contains $q$. So
$\{y\in[0,1)\colon y>q\}$ is the disjoint union of the intervals $J_k(q)$, $k\in\NN$. Hence $O_p$
is the disjoint union of the boxes $\prod_\ell J_{k_\ell}(p_\ell)$, $k\in\NN^d$; there is at most
one of each shape~$k$. Let $A_p$ be the union of those with $|k|>m^+$. The other boxes are
finitely many, and their indicators are fixed by $\Pi$, since they are constant on the boxes of
a shape $k$ with $|k|\le m^+$. So $\1_{O_p}-\Pi\1_{O_p}=\1_{A_p}-\Pi\1_{A_p}$.

For a shape $j$ with $|j|=m^+$ and a set $L$ of indices $\ell$ with $j_\ell\ge1$, let
$B_{j,L}$ be the product of the intervals $J_{j_\ell}(p_\ell)$ for $\ell\in L$ and of the dyadic
intervals of depth $j_\ell$ that contain $p_\ell$ for $\ell\notin L$; it is empty or a dyadic box
of shape $j$. Every box $\prod_\ell J_{k_\ell}(p_\ell)$ of $A_p$ lies in one of these boxes: choose
a shape $j\le k$ with $|j|=m^+$ and let $L$ be the set of the $\ell$ with $j_\ell=k_\ell$; for
$\ell\notin L$, the interval $J_{k_\ell}(p_\ell)$ lies in the dyadic interval of depth
$k_\ell-1\ge j_\ell$ that contains $p_\ell$. So $A_p$ is covered by at most
$2^d\binom{m^++d-1}{d-1}\le C(m+1)^{d-1}$ dyadic boxes, each of volume
$2^{-m^+}\le2^{-m}(m+2)^{-2d}$. Hence $|A_p|\le C(m+1)^{-d-1}2^{-m}$, and,
by Theorem~\ref{thm:W}(a) for the covering boxes, also $\omega(A_p)\le C(m+1)^{-d-1}2^{-m}$. By (a)
again, $\int g\,d\omega\le2\int g\,d\lambda$ for $g\ge0$ constant on the boxes of a shape $j$ with
$|j|\le m^+$. As $\|\Delta_j\1_{A_p}\|_1\le2^d|A_p|$ and there are $N_m\le C(m+1)^d$ shapes
$j$ with $|j|\le m^+$,
\begin{align*}
\Bigl|\int(\1_{A_p}-\Pi\1_{A_p})\,d\omega\Bigr|\ &\le\ \omega(A_p)+2\sum_{|j|\le m^+}\|\Delta_j\1_{A_p}\|_1\\
&\le\ \omega(A_p)+C(m+1)^d|A_p|\ \le\ \frac{C2^{-m}}{m+1}.
\end{align*}

\emph{The main part.} For $|j|\le m^+$, let $W_j=1+\sum_{s<j}f_s$, with $f_s$ as in
Theorem~\ref{thm:W}(b). We show that
\begin{equation}\label{eq:main}
\int\Pi V\,d\omega\ \ge\ \int V\,d\lambda+\sum_{|s|\le m^+-d}\int V\,f_s\,d\lambda-B.
\end{equation}
By definition, $\Pi V=\sum_{|j|\le m^+}\Delta_jV$. The function $\Delta_jV$ is constant on the boxes
of shape $j$ and bounded by $2^dB$, so Theorem~\ref{thm:W}(b) with $g=\Delta_jV$ gives
$|\int\Delta_jV\,d\omega-\int\Delta_jV\,W_j\,d\lambda|\le B/N_m$. Summing over the $N_m$ shapes $j$,
we get $\int\Pi V\,d\omega\ge\sum_j\int\Delta_jV\,W_j\,d\lambda-B$. By the symmetry of $\Delta_j$,
$\int\Delta_jV\,W_j\,d\lambda=\int V\,\Delta_jW_j\,d\lambda$. Since $f_s$ is a linear combination of
the $h_R$ with $r(R)=s$, the action of $\Delta_j$ on $1$ and on the $h_R$ gives $\Delta_0W_0=1$,
$\Delta_jW_j=f_{j-\1}$ if all entries of $j$ are positive, and $\Delta_jW_j=0$ otherwise. This
proves \eqref{eq:main}.

\emph{Conclusion.} Since $|V|\le B$, $\omega\ge0$ and $\omega(D)\le1$, the fine part and
\eqref{eq:main} give
\[
B\ \ge\ \int V\,d\omega\ \ge\ \int V\,d\lambda+\sum_{|s|\le m^+-d}\int V\,f_s\,d\lambda-B
-\frac{C2^{-m}}{m+1}\,|\nu|(D).
\]
We bound the terms on the right. First, $|\int V\,d\lambda|\le B$. The terms with $|s|<m$
vanish. By Lemma~\ref{lem:haar}(i) and the choice of $\sigma_R$, the terms with $|s|=m$ sum to
$\kappa\theta_m\sum_{R\in\cD_m}|\int\Lambda_R\,d\nu|$. For $|s|=m+k$ with $k\ge1$, the
coefficients of $f_s$ are bounded by $\|f_s\|_\infty$, since the $h_R$ with $r(R)=s$ have disjoint
supports and modulus one there. So, by Lemma~\ref{lem:haar} and Theorem~\ref{thm:W}(b),
\[
\Bigl|\int V\,f_s\,d\lambda\Bigr|\ \le\ \|f_s\|_\infty\sum_{r(R)=s}\Bigl|\int\Lambda_R\,d\nu\Bigr|
\ \le\ C(\kappa\theta_m)^2(k+1)^{2d}2^{-m-k}|\nu|(D).
\]
There are $\binom{m+k+d-1}{d-1}\le(m+1)^{d-1}(k+1)^{d-1}$ shapes with $|s|=m+k$, and
$\sum_k(k+1)^{3d-1}2^{-k}<\infty$. Collecting terms gives
\[
3B\ \ge\ \kappa\theta_m\sum_{R\in\cD_m}\Bigl|\int\Lambda_R\,d\nu\Bigr|
-C\Bigl(\kappa^2\theta_m^2(m+1)^{d-1}+\frac1{m+1}\Bigr)2^{-m}|\nu|(D),
\]
and we divide by~$3$.
\end{proof}

\section{Hard measures}\label{sec:hard}

In this section, we construct the hard measures $\mu_\eta$ (Lemma~\ref{lem:hard}) and separate
their distribution functions (Lemma~\ref{lem:sep}). For Lemma~\ref{lem:hard}, we need the following
small ball inequality for tent functions. It is \cite[Theorem~1.4]{W} for the function
$\varphi(y)=\min\{y,1-y\}$ on $[0,1)$; in the notation
there, $\varphi_R=|R|^{-1}\Lambda_R$. Theorem~\ref{thm:W} remains true if the constant $c$ in
\eqref{eq:theta} is decreased, since $\kappa$ can be multiplied by the same factor.

\begin{theorem}[Small ball inequality for tents]\label{thm:tent}
Let $d\ge3$, and let $\theta_m$ be as in \eqref{eq:theta}, with a possibly smaller constant
$c>0$. Then
\begin{equation}\label{eq:W1}
\Bigl\|\sum_{R\in\cD_m}a_R\Lambda_R\Bigr\|_\infty\ \ge\ \theta_m\,2^{-2m}\sum_{R\in\cD_m}|a_R|
\qquad\text{for all $m\in\NN_0$ and all real }a_R.
\end{equation}
\end{theorem}

\medskip\noindent\textbf{The grid.} Let $h_m=2^{-m-1}$, let
$\Gamma_m=\{0,h_m,\dots,1-h_m\}^d$ be the grid of mesh~$h_m$, and let
$\Gamma_m^\circ=\{h_m,2h_m,\dots,1-h_m\}^d$ be the set of its points without a zero coordinate.

\medskip\noindent\textbf{Point reflection.} Let $S(x)=\1-x$. For a dyadic interval
$I=[k2^{-j},(k+1)2^{-j})$, let $S(I)=[(2^j-k-1)2^{-j},(2^j-k)2^{-j})$, the mirror image of $I$ up to
its endpoints, and for a dyadic box let $S(R)=S(I_1)\times\dots\times S(I_d)$, a box of the same
shape. Then $\Lambda_I(1-y)=\Lambda_{S(I)}(y)$, since both vanish at the endpoints, and hence
$\Lambda_R\circ S=\Lambda_{S(R)}$. Since $S(I)=I$ only for $I=[0,1)$, no $R\in\cD_m$ with $m\ge1$
satisfies $S(R)=R$. The map $S$ is a bijection of $\Gamma_m^\circ$.

\medskip\noindent\textbf{Antisymmetric sign patterns.} We call
$\eta\in\{-1,1\}^{\cD_m}$ \emph{antisymmetric} if
$\eta_{S(R)}=-\eta_R$ for all $R\in\cD_m$. Such $\eta$ are determined by their values on a set
$\cD'_m$ that contains one box of each pair $\{R,S(R)\}$, so there are $2^{K_m/2}$ of them.

\medskip The next lemma constructs the hard measures. The separation theorem for convex
sets turns the inequality \eqref{eq:W1}, which is a statement about all tent sums, into a
statement about one measure.

\begin{lemma}[Hard measures]\label{lem:hard}
Let $d\ge3$, $m\in\NN$ and let $\eta\in\{-1,1\}^{\cD_m}$ be antisymmetric. There are a
probability measure $\mu_\eta$ with at most $K_m/2+1$ atoms, all in $\Gamma_m^\circ$, and a
number $c_\eta\ge2\theta_m2^{-2m}$ such that
\[
\int\Lambda_R\,d\mu_\eta-\int\Lambda_{S(R)}\,d\mu_\eta=c_\eta\,\eta_R\qquad\text{for all }R\in\cD_m.
\]
\end{lemma}
\begin{proof}
Let $F=\sum_{R\in\cD_m}a_R\Lambda_R$ with $a_R\in\RR$. The tent of an interval of length at least
$2^{-m}$ is affine on every interval $[c,c+h_m]$ with $c\in h_m\ZZ$. So on every cell
$c+[0,h_m]^d$, $c\in\Gamma_m$, the function $F$ is affine in each variable and attains its
maximum and its minimum over the cell at vertices. Every vertex lies in $\Gamma_m^\circ$ or
has a coordinate $0$ or $1$, where $F$ vanishes. We call $a$ antisymmetric if $a_{S(R)}=-a_R$ for
all $R$. Then $F\circ S=\sum_Ra_R\Lambda_{S(R)}=-F$, and since $S$ is a bijection of
$\Gamma_m^\circ$,
\begin{equation}\label{eq:maxF}
\|F\|_\infty=\max_{x\in\Gamma_m^\circ}F(x)\qquad\text{if $a$ is antisymmetric.}
\end{equation}

Let $A\subset\RR^{\cD_m}$ be the space of antisymmetric vectors; it has dimension $K_m/2$. For
$x\in\Gamma_m^\circ$, let $v(x)\in A$ have the entries
$v(x)_R=\Lambda_R(x)-\Lambda_{S(R)}(x)$, and let $C\subset A$ be the convex hull of the finite
set of these vectors. For $a\in A$, we have
$\langle a,v(x)\rangle=F(x)-F(S(x))=2F(x)$, so $\max_{y\in C}\langle a,y\rangle=2\|F\|_\infty$ by
\eqref{eq:maxF}. Let $w=2\theta_m2^{-2m}\eta\in A$. By $|\eta_R|=1$ and \eqref{eq:W1},
\[
\langle a,w\rangle\ \le\ 2\theta_m2^{-2m}\sum_R|a_R|\ \le\ 2\|F\|_\infty\ =\ \max_{y\in C}\langle a,y\rangle
\qquad\text{for all }a\in A.
\]
This gives $w\in C$, by the separation theorem for convex sets: otherwise, let $y_0$ be the point
of the compact convex set $C$ nearest to $w$, and let $a=w-y_0\ne0$. For $y\in C$ and $0<s\le1$, the point
$y_0+s(y-y_0)$ lies in $C$, so it is not closer to $w$ than $y_0$, that is,
\[
|a|^2\ \le\ |w-y_0-s(y-y_0)|^2\ =\ |a|^2-2s\langle a,y-y_0\rangle+s^2|y-y_0|^2.
\]
Dividing by $s$ and letting $s\to0$, we get $\langle a,y-y_0\rangle\le0$ for all $y\in C$. So
$\max_{y\in C}\langle a,y\rangle=\langle a,y_0\rangle=\langle a,w\rangle-|a|^2<\langle a,w\rangle$,
contrary to the display. By Carath\'eodory's theorem in $A$, $w=\sum_ic_iv(x_i)$ with $c_i\ge0$, $\sum_ic_i=1$ and at
most $K_m/2+1$ terms. The measure $\mu_\eta=\sum_ic_i\delta_{x_i}$ satisfies the claim with
$c_\eta=2\theta_m2^{-2m}$, since
$\int\Lambda_R\,d\mu_\eta-\int\Lambda_{S(R)}\,d\mu_\eta=\sum_ic_iv(x_i)_R=w_R=c_\eta\eta_R$.
\end{proof}

\begin{figure}[ht]\centering
\includegraphics[width=0.55\linewidth]{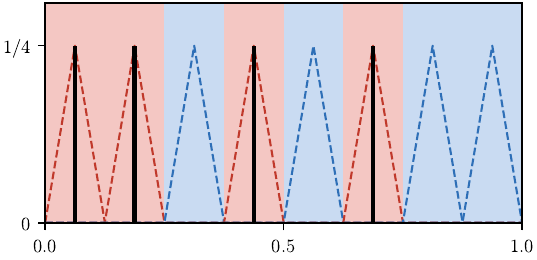}
\caption{$d=1$, $m=3$ and the antisymmetric $\eta=(1,1,-1,1,-1,1,-1,-1)$: the eight intervals
of $\cD_3$, red if $\eta_I=1$ and blue if $\eta_I=-1$, their tents (dashed, scaled), and
$\mu_\eta$, four atoms of mass $\frac14$. The reflection $x\mapsto1-x$ maps red intervals to blue
ones; $\mu_\eta$ has tent mass $\frac1{64}$ in each red interval and $0$ in each blue
one.}\label{fig:d1}
\end{figure}

\begin{figure}[t]\centering
\includegraphics[width=\linewidth]{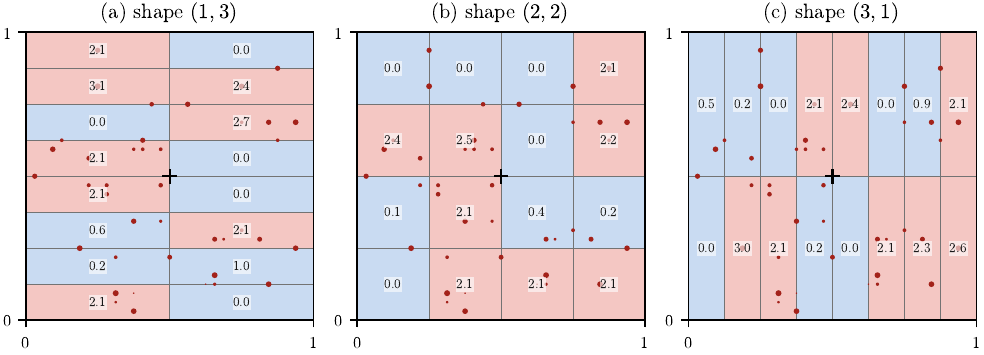}
\caption{$d=2$, $m=4$ and an antisymmetric $\eta$: boxes of three of the five shapes of
$\cD_4$, coloured by $\eta_R$ (red $+1$, blue $-1$), the measure $\mu_\eta$ (dots), and in each
box $R$ the tent mass of $\mu_\eta$ in $R$ divided by that of $\lambda$. The reflection through
the centre maps red boxes to blue ones; the red one has the larger tent mass, always by the same
amount.}\label{fig:lops}
\end{figure}

So the tent masses of $\mu_\eta$ in $R$ and in $S(R)$ differ by $c_\eta$, with the sign
$\eta_R$ (Figures~\ref{fig:d1} and~\ref{fig:lops}). For most $\eta$, the measure
$\mu_\eta$ is a hard measure for Theorem~\ref{thm:main}.

The reflection is needed for \eqref{eq:maxF}: for an antisymmetric tent sum $F$,
\eqref{eq:W1} is a lower bound for the maximum of $F$, not only of $|F|$. Without the
reflection, one would prescribe $\int\Lambda_R\,d\mu-\int\Lambda_R\,d\lambda=c\,\eta_R$ for all
$R$ and all $\eta$. This requires $\max F-\int F\,d\lambda\ge c\sum_R|a_R|$ for all tent sums
$F=\sum_Ra_R\Lambda_R$ (take $\eta_R=\sgn a_R$), which we do not have.

\medskip The next lemma separates the distribution functions of two hard measures.

\begin{lemma}[Separation]\label{lem:sep}
Let $d\ge3$. There are $c>0$ and $m_1\in\NN$ such that for $m\ge m_1$ and all antisymmetric
$\eta,\eta'\in\{-1,1\}^{\cD_m}$ that differ in at least $K_m/4$ entries, the measure
$\nu=\mu_\eta-\mu_{\eta'}$, with $\mu_\eta,\mu_{\eta'}$ as in Lemma~\ref{lem:hard}, satisfies
\begin{equation}\label{eq:sepinf}
\sup_{t\in D}|\nu([0,t])|\ \ge\ c\,(m+1)(\log(m+2))^{-2\zeta_d}\,2^{-m}.
\end{equation}
\end{lemma}
\begin{proof}
We have $|\nu|(D)\le2$. Let $x_R=\int\Lambda_R\,d\nu$. By Lemma~\ref{lem:hard},
$x_R-x_{S(R)}=c_\eta\eta_R-c_{\eta'}\eta'_R$. If $\eta_R\ne\eta'_R$, then $\eta'_R=-\eta_R$,
and so $|x_R-x_{S(R)}|=c_\eta+c_{\eta'}\ge4\theta_m2^{-2m}$. Since $S$ permutes $\cD_m$, and
$\eta_R\ne\eta'_R$ for at least $K_m/4$ boxes $R$,
\[
\sum_{R\in\cD_m}|x_R|\ =\ \frac12\sum_{R\in\cD_m}\bigl(|x_R|+|x_{S(R)}|\bigr)\ \ge\
\frac12\sum_{R\in\cD_m}|x_R-x_{S(R)}|\ \ge\ \frac{K_m}2\,\theta_m2^{-2m}.
\]
Let $C$ be the constant of Proposition~\ref{prop:sb}, let $c_1>0$ be such that
$K_m2^{-2m}\ge12c_1(m+1)^{d-1}2^{-m}$ for all $m$ (Section~\ref{sec:sb}), and let
$\kappa=\min\{1,c_1/(2C)\}$. By Proposition~\ref{prop:sb} and $|\nu|(D)\le2$,
\begin{align*}
\sup_{t\in D}|\nu([0,t])|\ &\ge\ \frac{\kappa\theta_m}3\cdot\frac{K_m}2\,\theta_m2^{-2m}
-2C\kappa^2\theta_m^2(m+1)^{d-1}2^{-m}-\frac{2C2^{-m}}{m+1}\\
&\ge\ c_1\kappa\,\theta_m^2(m+1)^{d-1}2^{-m}-\frac{2C2^{-m}}{m+1}.
\end{align*}
Here $\theta_m^2(m+1)^{d-1}$ is a constant times $(m+1)(\log(m+2))^{-2\zeta_d}$. So, for large $m$,
the first term is at least twice the second, and the claim follows.
\end{proof}

\begin{figure}[H]\centering
\includegraphics[width=0.7\linewidth]{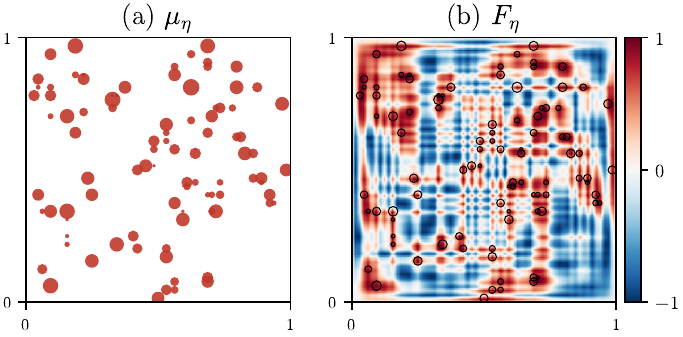}
\caption{$d=2$, $m=5$ and a random antisymmetric $\eta$. (a) $\mu_\eta$, computed by the linear program above,
with $95$ atoms; the area of a point is proportional to
its weight. (b) $F_\eta$ of Lemma~\ref{lem:loc}; it equals
$1$ at the atoms of $\mu_\eta$ (circled) and $-1$ at the mirror images of the atoms.}\label{fig:sigma}
\end{figure}

\medskip\noindent\textbf{Computing the hard measures.} The proof of Lemma~\ref{lem:hard}
gives a way to compute $\mu_\eta$. As there, let $A$ be the space of the $a\in\RR^{\cD_m}$ with
$a_{S(R)}=-a_R$ for all $R$, and for $x\in\Gamma_m^\circ$ let $v(x)\in A$ have the entries
$v(x)_R=\Lambda_R(x)-\Lambda_{S(R)}(x)$. Consider the linear program
\[
\text{maximise }c\quad\text{subject to}\quad c\,\eta=\sum_{x\in\Gamma_m^\circ}c_xv(x),\quad
c_x\ge0,\quad\sum_{x\in\Gamma_m^\circ}c_x=1,
\]
in the variables $c\in\RR$ and $c_x$, $x\in\Gamma_m^\circ$. By that proof, its value
satisfies $c^*\ge2\theta_m2^{-2m}$. By Carath\'eodory's theorem in $A$,
there is an optimal solution with at most $K_m/2+1$ nonzero $c_x$, and then
$\mu_\eta=\sum_xc_x\delta_x$ satisfies Lemma~\ref{lem:hard} with $c_\eta=c^*$. The dual program (see \cite[Section~6.2]{MG}) is
\[
\text{minimise }t\quad\text{subject to}\quad\langle a,\eta\rangle=1,\quad
\langle a,v(x)\rangle\le t\ \text{ for }x\in\Gamma_m^\circ,
\]
in the variables $t\in\RR$ and $a\in A$. As $\langle a,v(x)\rangle=2F(x)$ for
$F=\sum_Ra_R\Lambda_R$, by \eqref{eq:maxF} it minimises $2\|F\|_\infty$ over the $a\in A$ with
$\langle a,\eta\rangle=1$. The next lemma, which is not used in the proofs, shows that the atoms
of $\mu_\eta$ lie at maxima of $F$ for an optimal $a$; their mirror images lie at its minima, since
$F\circ S=-F$. Figure~\ref{fig:sigma} was computed in this way.

\begin{lemma}[Location of the atoms]\label{lem:loc}
Let $d\ge2$, $m\in\NN$ and let $\eta\in\{-1,1\}^{\cD_m}$ be antisymmetric. Let
$(c^*,(c_x))$ and $(t,a)$ be optimal solutions of the linear program above and of its dual
program, and let $F=\sum_Ra_R\Lambda_R$. Then $t=c^*=2\|F\|_\infty$, and
$F_\eta=F/\|F\|_\infty$ satisfies $F_\eta(x)=1$ whenever $c_x>0$.
\end{lemma}
\begin{proof}
For admissible $(c,(c_x))$ and $(t,a)$, since $\langle a,\eta\rangle=1$,
\[
c\ =\ \langle a,c\,\eta\rangle\ =\ \sum_xc_x\langle a,v(x)\rangle\ \le\ t\sum_xc_x\ =\ t.
\]
So the linear program is bounded, and by the duality theorem of linear programming \cite[Section~6.1]{MG}, the optimal
values agree: $c^*=t$. Hence, for the optimal solutions, equality holds above, and
$\langle a,v(x)\rangle=t$ whenever $c_x>0$. As $t$ is minimal, $t=\max_x\langle a,v(x)\rangle
=2\max_xF(x)=2\|F\|_\infty$ by \eqref{eq:maxF}. So $\|F\|_\infty=c^*/2$, which is positive, since $a\ne0$ and the tents
$\Lambda_R$, $R\in\cD_m$, are linearly independent, and $F(x)=\|F\|_\infty$ whenever $c_x>0$.
\end{proof}

\section{\texorpdfstring{Proofs of Theorem~\ref{thm:main} and Corollary~\ref{cor:one}}{Proofs of Theorem 1.1 and Corollary 1.2}}\label{sec:proof}

For a measure $\nu$, let $T\nu$ be its distribution function, $T\nu(t)=\nu([0,t])$ for
$t\in D$; for an integrable function $g$ we write $Tg$ for the distribution function of the
measure with the density~$g$. We identify a quadrature rule $Q_n$ with nodes $x_i\in D$ and
weights $w_i\in\RR$ with the measure $\sum_iw_i\delta_{x_i}$; then
$TQ_n(t)=\sum_{i\colon x_i\le t}w_i$. For a probability measure $\mu$ on $D$, the function
$T\mu-TQ_n$ is the \emph{discrepancy function} of $Q_n$ with respect to $\mu$. The
$L_p$-discrepancy of the introduction is its $L_p$-norm, $D_p(Q_n;\mu)=\|T\mu-TQ_n\|_p$ for
$1\le p<\infty$, and we also write $D_\infty(Q_n;\mu)=\|T\mu-TQ_n\|_\infty$. Then
$D_p(Q_n;\mu)\le D_\infty(Q_n;\mu)=D^*(Q_n;\mu)$, since distribution functions are
continuous from above. The smallest discrepancy with $n$ nodes is
\[
D_p(n;\mu)=\inf_{Q_n}D_p(Q_n;\mu)\qquad(1\le p\le\infty),
\]
the infimum over all rules with $n$ nodes and real weights. For a probability density $\rho$, we
write $D_p(Q_n;\rho)$ for the measure with this density.

The first lemma is the discretisation; it allows arbitrary weights. For $M\in\NN$ and $h=1/M$,
let $P_h$ be the mean value operator for the cells $c+[0,h)^d$, $c\in\{0,h,\dots,1-h\}^d$: the
function $P_hf$ is constant on each cell, and equal to the mean value of $f$ there. It is linear and $\|P_hf\|_p\le\|f\|_p$.
For $t$ in the cell $c+[0,h)^d$, let $U_t$ be the union of the $d$ slabs
$\{x\in D\colon c_\ell\le x_\ell<c_\ell+h\}$, each of volume $h$. If $t'$ lies in the same cell,
then $[0,t]$ and $[0,t']$ differ by a subset of $U_t$.

\begin{lemma}[Few step functions]\label{lem:few}
Let $n\ge2$, $M\in\NN$, $h=1/M$, $1\le p\le\infty$ and $n^{-2}\le\tau\le2$. Let $\cG$ be a set of
probability measures on $D$ with $\|T\mu-P_hT\mu\|_p\le\frac12n^{-2}$ for all $\mu\in\cG$. There
is a set $\cC$ of at most $M^{dn}(1+6n^2)^{2^dn}$ bounded functions on $D$ such that for every
$\mu\in\cG$ with $D_p(n;\mu)<\tau$ there is $y\in\cC$ with $\|T\mu-y\|_p<3\tau$.
\end{lemma}

The assumption holds if every $\mu\in\cG$ is supported on the grid $\{0,h,\dots,1-h\}^d$, since
then $T\mu$ is constant on the cells. It also holds if $\mu(U_t)\le\frac12n^{-2}$ for all $\mu\in\cG$
and $t\in D$; in particular if every $\mu\in\cG$ has a density
bounded by $B$ and $M\ge2dBn^2$.

\begin{proof}
The distribution function of a rule is $TQ_n=\sum_iw_i\1_{[x_i,1)}$. In one variable, if
$x=c+\vartheta h$ with $c\in\{0,h,\dots,1-h\}$ and $0\le\vartheta<1$, then
$P_h\1_{[x,1)}=(1-\vartheta)\1_{[c,1)}+\vartheta\1_{[c+h,1)}$: both sides vanish to the left
of the cell $[c,c+h)$, equal $1-\vartheta$ on it, and are one to the right of it. Taking
products, $P_h\1_{[x,1)}$ is a linear combination of the $2^d$ functions $\1_{[z,1)}$,
$z\in c+\{0,h\}^d$, if $x$ lies in the cell with corner $c$. Each of the $n$ nodes lies in
one of the $M^d$ cells, so there are at most $M^{dn}$ possibilities for the cells
$c_i+[0,h)^d$ that contain the nodes $x_i$, $i=1,\dots,n$. For each of them, let
$Z=\bigcup_i(c_i+\{0,h\}^d)$ be the set of the at most $2^dn$ vertices of these cells, and let
$V_Z$ be the span of the functions $\1_{[z,1)}$, $z\in Z$. Then
$\dim V_Z\le2^dn$, and $P_hTQ_n\in V_Z$ for every rule with nodes in these cells, whatever the
weights are. In each $V_Z$, choose a maximal subset of $\{g\in V_Z\colon\|g\|_p\le3\}$ with mutual
distances at least $n^{-2}$. The balls of radius $\frac12n^{-2}$ in $V_Z$ around its elements are
disjoint and lie in the ball of radius $3+\frac12n^{-2}$. Comparing the volumes of these balls
in $V_Z$, the subset has at most
\[
\Bigl(\frac{3+\frac12n^{-2}}{\frac12n^{-2}}\Bigr)^{\dim V_Z}\ \le\ (1+6n^2)^{2^dn}
\]
elements. Let $\cC$ be the union of these subsets over all possibilities for the cells.

Let $\mu\in\cG$ and let $Q_n$ be a rule with $\|T\mu-TQ_n\|_p<\tau$. Then
$\|P_hTQ_n\|_p<\|P_hT\mu\|_p+\tau\le3$, so, by maximality, some $g\in\cC$ has
$\|P_hTQ_n-g\|_p<n^{-2}$, and
\[
\|T\mu-g\|_p\le\|T\mu-P_hT\mu\|_p+\|P_h(T\mu-TQ_n)\|_p+\|P_hTQ_n-g\|_p
<\tau+\tfrac32n^{-2}\ \le\ \tfrac52\tau.\qedhere
\]
\end{proof}

The second lemma is the counting. The measure $\nu$ and the weight $\xi$ in it are needed only
for Corollary~\ref{cor:one}; in the main proof $\xi=1$ and $\nu=0$.

\begin{lemma}[Most sign patterns are hard]\label{lem:count}
In the situation of Lemma~\ref{lem:few}, let $m\in\NN$ and $0<\xi\le1$, let $\nu\ge0$ be a
measure on $D$ with $\nu(D)=1-\xi$, and let $\mu_\eta$, for antisymmetric
$\eta\in\{-1,1\}^{\cD_m}$, be probability measures on $D$ such that the probability measures
$\pi_\eta=\nu+\xi\mu_\eta$ lie in $\cG$. Assume that $\log_2\#\cC\le K_m/24$ and that
\[
\xi\,\|T\mu_\eta-T\mu_{\eta'}\|_p\ \ge\ 6\tau
\]
whenever $\eta,\eta'$ differ in at least $K_m/4$ entries. Then $D_p(n;\pi_\eta)\ge\tau$ for all antisymmetric $\eta$
outside a set of proportion at most $2^{-K_m/20}$.
\end{lemma}
\begin{proof}
Recall from Section~\ref{sec:hard} that $\cD'_m$ contains one box of each pair
$\{R,S(R)\}$, so it has $K=K_m/2$ elements. We identify an antisymmetric $\eta$ with its
values on $\cD'_m$, an element of $\{-1,1\}^{\cD'_m}$; since $\eta_{S(R)}=-\eta_R$, two
antisymmetric $\eta,\eta'$ differ in exactly twice as many boxes of $\cD_m$ as of $\cD'_m$. Let
$E$ be the set of antisymmetric $\eta$ with $D_p(n;\pi_\eta)<\tau$. By Lemma~\ref{lem:few}, each
$\eta\in E$ has some $y_\eta\in\cC$ with $\|T\pi_\eta-y_\eta\|_p<3\tau$. If
$y_\eta=y_{\eta'}$, then, as $T\pi_\eta-T\pi_{\eta'}=\xi(T\mu_\eta-T\mu_{\eta'})$,
\[
\xi\|T\mu_\eta-T\mu_{\eta'}\|_p<6\tau,
\]
so, by the assumption, $\eta$ and $\eta'$ differ in fewer than $K_m/4$ boxes of $\cD_m$, hence
in fewer than $K/4$ boxes of $\cD'_m$. So, for each $y\in\cC$, the $\eta\in E$ with $y_\eta=y$
lie in the Hamming ball of radius $K/4$ in $\{-1,1\}^{\cD'_m}$ around any of them. Such a ball has $N=\sum_{i\le K/4}\binom Ki$ elements, and
\[
1\ \ge\ \sum_{i\le K/4}\binom Ki\bigl(\tfrac14\bigr)^i\bigl(\tfrac34\bigr)^{K-i}\ \ge\
N\bigl(\tfrac14\bigr)^{K/4}\bigl(\tfrac34\bigr)^{3K/4}
\]
gives $N\le2^{0.812K}$. Therefore $\#E\le\#\cC\,2^{0.812K}\le2^{0.9K}$, since
$\log_2\#\cC\le K/12$. This is the proportion $2^{-0.1K}=2^{-K_m/20}$ of the $2^K$ antisymmetric
sign patterns.
\end{proof}

\begin{proof}[\normalfont\bfseries Proof of Theorem~\ref{thm:main}]
Let $c$ and $m_1$ be as in Lemma~\ref{lem:sep}, and let
\[
s_m=c\,(m+1)(\log(m+2))^{-2\zeta_d}\,2^{-m}
\]
be the right-hand side of \eqref{eq:sepinf}.

\emph{Parameters.} Let $C_0=2d+7\cdot2^d$ and $b_d=1/(24\,C_0\,(d-1)!)$.
Given a large $n$,
let $m$ be the least integer with $b_d2^mm^{d-2}\ge n$. Then $m\ge m_1$,
$\frac12\log_2n\le m\le C\log n$ and $b_d2^{m-1}(m-1)^{d-2}<n\le b_d2^mm^{d-2}$. In particular,
$2^{-m}>b_d(m-1)^{d-2}/(2n)$, and inserting this into $s_m$ gives
\begin{equation}\label{eq:sm}
s_m\ \ge\ c'\,\frac{(\log n)^{d-1}}{n(\log\log n)^{\chi_d}},\qquad\chi_d=2\zeta_d .
\end{equation}
Let $M=2^{m+1}$. Since $\log_2M\le2m$ and $\log_2(1+6n^2)\le7m$, the set $\cC$ of
Lemma~\ref{lem:few} satisfies
\begin{equation}\label{eq:C0}
\log_2\#\cC\le dn\log_2M+2^dn\log_2(1+6n^2)\le C_0\,n\,m\le C_0b_d2^mm^{d-1}\le K_m/24.
\end{equation}

\emph{Counting.} Let $p=\infty$, $h=1/M=h_m$, let $\cG$ be the set of probability measures supported on
$\Gamma_m$, and $\tau=s_m/6$. By the remark after Lemma~\ref{lem:few}, $\cG$ satisfies the
assumption of that lemma. For antisymmetric $\eta$, the measures $\mu_\eta$ of
Lemma~\ref{lem:hard} lie in $\cG$. Distribution functions of measures in $\cG$ are constant on
the cells $c+[0,h)^d$, $c\in\Gamma_m$, so their $L_\infty$-norm is their supremum, and by
Lemma~\ref{lem:sep}, $\|T\mu_\eta-T\mu_{\eta'}\|_\infty\ge s_m$ if $\eta,\eta'$ differ in at
least $K_m/4$ entries.
By \eqref{eq:C0}, $\log_2\#\cC\le K_m/24$, and $n^{-2}\le\tau\le2$ for large $n$
by \eqref{eq:sm}, since $s_m\to0$. Lemma~\ref{lem:count} with $\xi=1$ and
$\nu=0$ gives $D_\infty(n;\mu_\eta)\ge\tau$, that is, $D^*(Q_n;\mu_\eta)\ge\tau$ for every rule
$Q_n$ with $n$ nodes, for all antisymmetric $\eta$ outside a set of proportion $2^{-K_m/20}$.
With \eqref{eq:sm} this proves the theorem for $\mu=\mu_\eta$. This measure has at most
$K_m/2+1$ atoms, and by the choice of $m$,
$K_m\le2^m(m+1)^{d-1}<2n(m+1)^{d-1}/(b_d(m-1)^{d-2})\le C_dn\log n$.
\end{proof}

The next proposition is the first part of the proof of Corollary~\ref{cor:one}.

\begin{proposition}[One measure for infinitely many $n$]\label{prop:one}
Let $d\ge3$ and let $(\Phi_n)$ be a sequence of positive numbers with $\Phi_n\to\infty$. There
are $c_d>0$, not depending on $(\Phi_n)$, and a probability measure $\mu$ on $[0,1)^d$ with
countable support such that, for infinitely many $n$, every quadrature rule $Q_n$ with $n$
nodes in $[0,1)^d$ and real weights satisfies
\[
D^*(Q_n;\mu)\ \ge\ c_d\,\frac{(\log n)^{d-1}}{n\,\Phi_n\,(\log\log n)^{\chi_d}} .
\]
\end{proposition}
\begin{proof}
We use the proof of Theorem~\ref{thm:main}. For large $m$, let $N(m)=\lfloor b_d2^mm^{d-2}\rfloor$. Then $m$ is the least integer with
$b_d2^mm^{d-2}\ge N(m)$, since $b_d2^{m-1}(m-1)^{d-2}<b_d2^mm^{d-2}-1$, so the preceding proof applies to $n=N(m)$ and this $m$; in
particular, by \eqref{eq:sm}, $s_m\ge1/N(m)$.

Choose $m(1)<m(2)<\dots$ so large that this holds and that $N_j=N(m(j))$ satisfies
$N_j\ge64N_{j-1}^3$ and $\Phi_{N_j}\ge N_{j-1}^3$ for $j\ge2$. Let $\xi_j=N_{j-1}^{-3}$ for
$j\ge2$, $\xi_1=1-\sum_{j\ge2}\xi_j>0$ and $\bar\xi_j=\sum_{l>j}\xi_l\le2N_j^{-3}$. Then
$\Gamma_{m(l)}\subset\Gamma_{m(j)}$ for $l<j$. We choose probability measures $\mu_1,\mu_2,\dots$,
where $\mu_j$ is supported on $\Gamma_{m(j)}$, one after the other. Let $\mu_1=\delta_0$. For $j\ge2$, let
$\nu_j=\sum_{l<j}\xi_l\mu_l+\bar\xi_j\delta_0$; then $\nu_j\ge0$, $\nu_j(D)=1-\xi_j$, and
$\nu_j$ is supported on $\Gamma_{m(j)}$. We apply Lemma~\ref{lem:count} with $p=\infty$, $n=N_j$,
$m=m(j)$, the measures $\mu_\eta$ for antisymmetric $\eta$, $\xi=\xi_j$, $\nu=\nu_j$, $\tau=\xi_js_m/6$, and with $\cG$ and $M$ as in
the preceding proof. The measures $\nu_j+\xi_j\mu_\eta$ are supported on $\Gamma_{m(j)}$, so they
lie in $\cG$; the bound $\log_2\#\cC\le K_m/24$ holds as before; and $\tau\ge n^{-2}$,
since $6N_j^{-2}<N_{j-1}^{-3}N_j^{-1}\le\xi_js_m$, and $\tau\le s_m/6\le2$. So there is an antisymmetric $\eta$ such that
$\mu_j:=\mu_\eta$ satisfies $D_\infty(n;\nu_j+\xi_j\mu_j)\ge\xi_js_m/6$.

Let $\mu=\sum_j\xi_j\mu_j$, a probability measure with countable support. For $j\ge2$,
$\mu-(\nu_j+\xi_j\mu_j)=\sum_{l>j}\xi_l\mu_l-\bar\xi_j\delta_0$ is the difference of two
nonnegative measures of total mass $\bar\xi_j$, so
$|T\mu(t)-T(\nu_j+\xi_j\mu_j)(t)|\le\bar\xi_j\le2N_j^{-3}\le\xi_js_m/12$ for all $t$. Hence,
for $n=N_j$, $m=m(j)$ and every rule $Q_n$ with $n$ nodes and real weights, by \eqref{eq:sm},
\[
D^*(Q_n;\mu)\ \ge\ \frac{\xi_js_m}{12}\ \ge\ \frac{s_m}{12\,\Phi_n}\ \ge\
c_d\,\frac{(\log n)^{d-1}}{n\,\Phi_n(\log\log n)^{\chi_d}} .
\]
\end{proof}

\begin{proof}[\normalfont\bfseries Proof of Corollary~\ref{cor:one}]
Let $C$ be the constant of \eqref{eq:Wup} and apply Proposition~\ref{prop:one} with
\[
\Phi_n=(c_d/C)\,(\log n)^{\delta}(\log\log n)^{-\chi_d}\qquad(n\ge3);
\]
this is allowed since $c_d$ does not depend on $(\Phi_n)$. Then the right-hand side in
Proposition~\ref{prop:one} is
\[
C\,\frac{(\log n)^{d-1-\delta}}n\ \ge\ (\log n)^{\gamma_d-\delta}\inf_{\cP_n}D^*(\cP_n;\lambda),
\]
where the inequality holds by \eqref{eq:Wup}.
\end{proof}

\section{\texorpdfstring{The $L_p$-discrepancy}{The Lp-discrepancy}}\label{sec:L2}

For the $L_p$-discrepancy, $1<p<\infty$, the small ball inequality for all coefficients is
not needed. The bounds for entropy numbers of Bilyk, Neumayer and T.~Ullrich \cite{BNU}, which rest
on a small ball inequality for Gaussian coefficients, and the duality of metric entropy of
Artstein, Milman and Szarek \cite{AMS} and of Artstein, Milman, Szarek and
Tomczak-Jaegermann \cite{AMST} suffice.
Section~\ref{sec:L2thm} treats the worst measure and the star discrepancy for densities in
$L_r$, Section~\ref{sec:L2expl} an elementary lower bound for densities that are bounded in $L_2$, and Section~\ref{sec:L2dens} an upper bound for such densities.

\subsection{The worst measure}\label{sec:L2thm}

The following theorem is the lower bound of Theorem~\ref{thm:L2intro}, with a measure that
has a density bounded by a power of $n$.

\begin{theorem}[$L_p$, worst measure]\label{thm:L2meas}
Let $d\ge2$ and $1<p<\infty$. There are $c>0$, $C<\infty$ and $n_0\in\NN$ such that for every
$n\ge n_0$ there is a probability density $\rho\le n^C$ on $D$ such that every quadrature
rule $Q_n$ with $n$ nodes and real weights satisfies
\[
D_p(Q_n;\rho)\ \ge\ c\,\frac{(\log n)^{d-3/2}}{n} .
\]
\end{theorem}

The proof uses entropy numbers. For a subset $F$ of a normed space $Y$ and a real number
$k\ge1$, the entropy number $e_k(F,Y)$ is the infimum of all $r>0$ such that $F$ is covered by at
most $2^{k-1}$ balls of radius $r$ in $Y$; see Carl and Stephani \cite{CS}. For
$1\le r\le\infty$, let $B_r$ be the unit ball of $L_r$, and for $1<p<\infty$ let $p'=p/(p-1)$. We consider the
operator $T\colon L_1\to L_p$, $g\mapsto Tg$, with $\|Tg\|_{p}\le\|g\|_1$, and its dual operator
$T^*\colon L_{p'}\to L_\infty$, $T^*f(x)=\int_{[x,1)}f\,d\lambda$, and let
$e_k(T)=e_k(T(B_1),L_p)$ and $e_k(T^*)=e_k(T^*(B_{p'}),L_\infty)$. For $1<p\le2$, there
are $a,b\ge1$, depending only on $p$, such that
\begin{equation}\label{eq:AMS}
e_{bk}(T^*)\ \le\ a\,e_k(T)\qquad\text{for all }k\ge1.
\end{equation}
For $p=2$, this is the duality of metric entropy of Artstein, Milman and Szarek
\cite[Theorem~1]{AMS}, with universal constants; the operator form is explained in the
introduction and after Observation~4 there. For $1<p<2$, it follows in the same way from the
duality of metric entropy for $K$-convex spaces such as $L_p$ of Artstein, Milman, Szarek and
Tomczak-Jaegermann \cite[Theorem~5]{AMST}; see also Szarek \cite[Theorem~2.2]{Sz}. By duality, the bound of Bilyk, Neumayer
and T.~Ullrich \cite[Theorem~1.2]{BNU} for the entropy numbers of $T\colon L_{p'}\to
C([0,1]^d)$ gives the following lower bound.

\begin{lemma}[Entropy of distribution functions]\label{lem:entropy}
Let $d\ge2$ and $1<p<\infty$. There is $c>0$ such that for all real $k\ge2$
\[
e_k\bigl(T(B_1),L_p\bigr)\ \ge\ c\,\frac{(\log k)^{d-1/2}}k .
\]
\end{lemma}
\begin{proof}
For $p>2$, the claim follows from the case $p=2$, since $\|\cdot\|_p\ge\|\cdot\|_2$. Let
$1<p\le2$. Recall that $Tg(t)=\int_{[0,t]}g\,d\lambda$. With the point reflection $S(x)=\1-x$, we have
$T^*f=(T(f\circ S))\circ S$, and composition with $S$ is an isometry of $L_{p'}$ and of
$L_\infty$ that maps $C([0,1]^d)$ onto itself. So $T(B_{p'})=\{h\circ S\colon h\in
T^*(B_{p'})\}$. Let $\varrho>e_k(T)$. By \eqref{eq:AMS}, $e_{bk}(T^*)<a\varrho$, so
$T^*(B_{p'})$ is covered by at most $2^{bk-1}$ balls of radius $a\varrho$ in $L_\infty$. We
replace each of these balls that meets $T^*(B_{p'})$ by the ball of radius $2a\varrho$ around one
of its points in $T^*(B_{p'})\subset C([0,1]^d)$, and compose with $S$. So $T(B_{p'})$ is
covered by at most $2^{j-1}$ balls of radius $2a\varrho$ with centres in $C([0,1]^d)$, where
$j=\lceil bk\rceil$; that is, $e_j(T(B_{p'}),C([0,1]^d))\le2a\varrho$. By
\cite[Theorem~1.2]{BNU} with the exponent $p'$,
$e_j(T(B_{p'}),C([0,1]^d))\ge c_0\,j^{-1}(\log j)^{d-1/2}$ for all integers $j\ge2$.
Since $k\le j\le2bk$, we get $2a\varrho\ge c_0(2bk)^{-1}(\log k)^{d-1/2}$. Letting
$\varrho\to e_k(T)$ gives the claim with $c=c_0/(4ab)$.
\end{proof}

The next lemma bounds the worst-case discrepancy from below by the entropy number $e_k$ with
$k\approx n\log n$.

\begin{lemma}[Discrepancy and entropy]\label{lem:discent}
Let $d\ge2$, $1\le p\le\infty$ and $1\le r\le\infty$, with $r>1$ if $p=\infty$. There is
$A>0$, depending only on $d$, $p$ and $r$, such that for all $n\ge2$
\[
\sup_{\mu\in\cG_r}D_p(n;\mu)\ \ge\ \frac17\,e_k\bigl(T(B_r),L_p\bigr)-n^{-2},
\qquad k=\lceil An\log_2n\rceil,
\]
where $\cG_r$ is the set of probability measures on $D$ with a density $\rho$ with
$\|\rho\|_r\le2$.
\end{lemma}

We use the lemma for $r=1$ and
$1<p<\infty$, where Lemma~\ref{lem:entropy} bounds the entropy numbers, and for $1<r<\infty$ and
$p=\infty$.

\begin{proof}
Let $M=\lceil d(2n^2)^p\rceil$ if $p<\infty$, $M=\lceil d(4n^2)^{r/(r-1)}\rceil$ if $p=\infty$,
and $h=1/M$. First, $\cG_r$ satisfies the assumption of Lemma~\ref{lem:few}. Let $\mu\in\cG_r$
have the density $\rho$. With $U_t$ as in
Section~\ref{sec:proof},
$|T\mu-P_hT\mu|(t)\le\mu(U_t)\le1$. If $p<\infty$: a point $x$ lies in $U_t$ only for
$t$ in the $d$ slabs through the cell of $x$, so, by Fubini's theorem,
\[
\|T\mu-P_hT\mu\|_p^p\ \le\ \|T\mu-P_hT\mu\|_1\ \le\ \int_D\mu(U_t)\,dt\ \le\ dh\ \le\ (2n^2)^{-p}.
\]
If $p=\infty$: H\"older's inequality gives $\mu(U_t)\le\|\rho\|_r(dh)^{1-1/r}\le\frac12n^{-2}$.

Let $\sigma=\sup_{\mu\in\cG_r}D_p(n;\mu)$; then $\sigma\le1$, since $\|T\mu\|_p\le1$. Let
$\max\{\sigma,n^{-2}\}<\tau\le2$, let $\cC$ be as in Lemma~\ref{lem:few} for $\cG=\cG_r$, and let $\cC_1$ be the
set of all functions $y-y'+c'\,T\lambda$ with $y,y'\in\cC$, $c'\in n^{-2}\ZZ$ and $|c'|\le1$. Let
$g\in B_r$. Then $\|g\|_1\le\|g\|_r\le1$, so $0\le\int g_\pm\,d\lambda\le1$, and the
measures $\mu_\pm=g_\pm\,d\lambda+(1-\int g_\pm\,d\lambda)\lambda$ are probability measures with the
densities $g_\pm+1-\int g_\pm\,d\lambda$. Their $L_r$-norms are at most $\|g_\pm\|_r+1\le2$, so
they lie in $\cG_r$. Moreover,
$Tg=T\mu_+-T\mu_-+(\int g\,d\lambda)\,T\lambda$. Lemma~\ref{lem:few} gives $y,y'\in\cC$ with
$\|T\mu_+-y\|_p<3\tau$ and $\|T\mu_--y'\|_p<3\tau$. With the nearest point $c'$
of $n^{-2}\ZZ$ to $\int g\,d\lambda\in[-1,1]$, and as $\|T\lambda\|_p\le1$, we get
$\|Tg-(y-y'+c'T\lambda)\|_p<7\tau$. So the balls of radius $7\tau$ around the
elements of $\cC_1$ cover $T(B_r)$. By Lemma~\ref{lem:few},
$\log_2\#\cC\le dn\log_2M+2^dn\log_2(1+6n^2)$, so
$\log_2\#\cC_1\le2\log_2\#\cC+\log_2(2n^2+1)\le k-1$ for $k=\lceil An\log_2n\rceil$ with a
suitable $A$. Hence $e_k(T(B_r),L_p)\le7\tau$, and letting $\tau\to\max\{\sigma,n^{-2}\}$ gives
$e_k(T(B_r),L_p)\le7\max\{\sigma,n^{-2}\}\le7\sigma+7n^{-2}$, which is the claim.
\end{proof}

\begin{proof}[\normalfont\bfseries Proof of Theorem~\ref{thm:L2meas}]
Let $k=\lceil An\log_2n\rceil$, with $A$ as in Lemma~\ref{lem:discent} for $r=1$. By
Lemma~\ref{lem:entropy}, $e_k(T(B_1),L_p)\ge c'k^{-1}(\log k)^{d-1/2}$. So Lemma~\ref{lem:discent}
gives
\[
\sup_{\mu\in\cG_1}D_p(n;\mu)\ \ge\ \frac{c'}7\,\frac{(\log k)^{d-1/2}}k-n^{-2}\ \ge\
2c\,\frac{(\log n)^{d-3/2}}n
\]
for large $n$. So some probability density $\rho$ has
$D_p(n;\rho)\ge{}\frac32c(\log n)^{d-3/2}/n$. Let $M=\lceil2dn^{2p}\rceil$, $h=1/M$, and let
$\tilde\rho=P_h\rho$ be the mean value of $\rho$ on each cell $c+[0,h)^d$, with $P_h$ as in
Section~\ref{sec:proof}. Then $\tilde\rho$ is a probability density and $\tilde\rho\le M^d\le n^C$
for large $n$, with $C$ depending only on $d$ and $p$. For $t$ in a cell, the cells that lie in $[0,t]$
or do not meet it contribute the same to $T\rho(t)$ and $T\tilde\rho(t)$, and the other cells lie
in $U_t$ (Section~\ref{sec:proof}). So
$|T\rho(t)-T\tilde\rho(t)|\le2\int_{U_t}\rho\,d\lambda$. As in the proof of
Lemma~\ref{lem:discent}, $\int_D\int_{U_t}\rho\,d\lambda\,dt\le dh$. As $|T\rho-T\tilde\rho|\le1$, we get
\[
\|T\rho-T\tilde\rho\|_p^p\ \le\ \|T\rho-T\tilde\rho\|_1\ \le\ 2dh\ \le\ n^{-2p}.
\]
So $D_p(n;\tilde\rho)\ge D_p(n;\rho)-n^{-2}\ge c(\log n)^{d-3/2}/n$ for large $n$.
\end{proof}

\begin{proof}[\normalfont\bfseries Proof of Theorem~\ref{thm:L2intro}]
The lower bound is Theorem~\ref{thm:L2meas}; it holds for all weights, in particular for equal
weights. The upper bound follows from \cite[Theorem~1]{ABN}, since the $L_p$-norm is at most the
supremum.
\end{proof}

The $L_r$-norms of the densities of Theorem~\ref{thm:L2meas}, $r>1$, may grow with $n$. The
next proposition gives the exponent $d-\frac32$ for the star discrepancy and densities with
$\|\rho\|_r\le2$; for $r=2$ this is an entry of Table~\ref{tab}. It follows from
Lemma~\ref{lem:discent} with $p=\infty$ and the bound of Bilyk, Neumayer and T.~Ullrich
\cite[Theorem~1.2]{BNU} for $T\colon L_r\to C([0,1]^d)$, without duality.

\begin{proposition}[Star discrepancy for densities in $L_r$]\label{prop:Lr}
Let $d\ge2$ and $1<r<\infty$. There are $c>0$ and $n_0\in\NN$ such that for every $n\ge n_0$
there is a probability density $\rho$ on $D$ with $\|\rho\|_r\le2$ such that every quadrature rule
$Q_n$ with $n$ nodes and real weights satisfies
\[
D_\infty(Q_n;\rho)\ \ge\ c\,\frac{(\log n)^{d-3/2}}n .
\]
\end{proposition}
\begin{proof}
The functions $Tg$, $g\in B_r$, are continuous on $[0,1]^d$, so the $L_\infty$-norm of their
differences is the supremum on $[0,1]^d$. If $T(B_r)$ is covered by balls in $L_\infty$, we
replace each ball that meets $T(B_r)$ by the ball of twice the radius around one of its points in
$T(B_r)$; so $e_k(T(B_r),C([0,1]^d))\le2e_k(T(B_r),L_\infty)$. By \cite[Theorem~1.2]{BNU} for
$p=r$, the left-hand side is at least $c_0k^{-1}(\log k)^{d-1/2}$ for large $k$. With
Lemma~\ref{lem:discent} for $p=\infty$ and $k=\lceil An\log_2n\rceil$, we get
\[
\sup_{\mu\in\cG_r}D_\infty(n;\mu)\ \ge\ \frac{c_0}{14}\,\frac{(\log k)^{d-1/2}}k-n^{-2}\ \ge\
2c\,\frac{(\log n)^{d-3/2}}n
\]
for large $n$, and the claim follows.
\end{proof}

\subsection{\texorpdfstring{An elementary bound for densities in $L_2$}{An elementary bound for densities in L2}}\label{sec:L2expl}

\begin{figure}[ht]\centering
\includegraphics[width=0.9\linewidth]{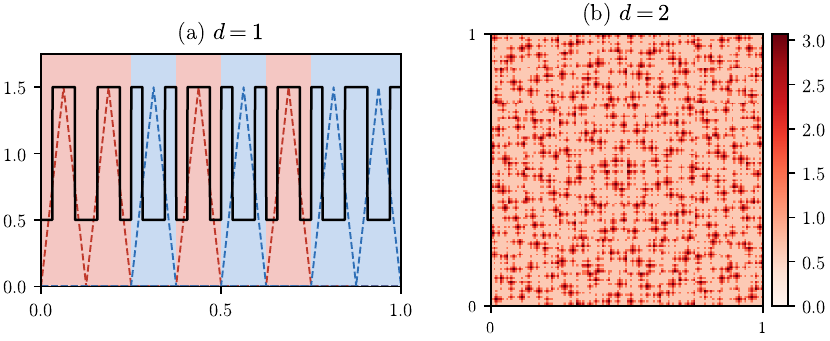}
\caption{The densities $\rho_\eta$. (a) $d=1$, $m=3$ and $\eta$ as in Figure~\ref{fig:d1}:
$\rho_\eta$ is large where the tents of the red intervals are high and where those of the blue
intervals are low; its tent mass is $1.25$ times that of $\lambda$ in the red intervals and
$0.75$ times in the blue ones. (b) $d=2$, $m=5$ and the $\eta$ of
Figure~\ref{fig:sigma}.}\label{fig:g}
\end{figure}

In this section, we give densities $\rho_\eta$ with $\|\rho_\eta\|_2\le2$ such that,
for most $\eta$,
every rule with $n$ nodes and real weights has $L_2$-discrepancy at least
$c(\log n)^{d-2}/n$ with respect to $\rho_\eta$. The proof is elementary: it uses neither the results of
Wengiel nor those of Bilyk, Neumayer and T.~Ullrich. As for $\mu_\eta$, their tent masses in $R$ and $S(R)$ differ with the sign
$\eta_R$, but only by $8^{-d}(m+1)^{-(d-1)/2}|R|^2$ instead of $2\theta_m|R|^2$
(Figure~\ref{fig:g}).

Let $u$ take the values $-1,1,1,-1$ on the four quarters of $[0,1)$ and $0$
elsewhere, let $\psi_I(y)=u((y-\inf I)/|I|)$ for a dyadic interval $I$ and
$\psi_R(x)=\prod_{\ell}\psi_{I_{\ell}}(x_{\ell})$ for a dyadic box $R=I_1\times\dots\times I_d$, and for
$\eta\in\{-1,1\}^{\cD_m}$ let
\[
g_\eta=(m+1)^{-\frac{d-1}2}\sum_{R\in\cD_m}\eta_R\psi_R,\qquad
\rho_\eta=1+(g_\eta)_+-\int_D(g_\eta)_+\,d\lambda .
\]
For antisymmetric $\eta$, we have $g_\eta\circ S=-g_\eta$ and
$\rho_\eta-\rho_\eta\circ S=g_\eta$.

\begin{proposition}[Elementary lower bound for densities in $L_2$]\label{prop:dens}
Let $d\ge2$. There are $b_d',c_d>0$ and $n_d\in\NN$ with the following property. Let $n\ge n_d$
and let $m$ be the least integer with $b_d'2^mm^{d-2}\ge n$. Then all $\rho_\eta$ are probability
densities with $\|\rho_\eta\|_2\le2$, and for all antisymmetric $\eta$ outside a set of
proportion $2^{-K_m/20}\le2^{-c_dn\log n}$, every rule $Q_n$ with $n$ nodes and real
weights satisfies
\[
D_2(Q_n;\rho_\eta)\ \ge\ c_d\,\frac{(\log n)^{d-2}}n .
\]
\end{proposition}
\begin{proof}
Let $b_d'=2^{-d-5}/(24\,(d-1)!)$.
In one variable, $\int_0^1\min\{y,1-y\}u(y)\,dy=\frac18$, and $\psi_I$ is orthogonal to the
affine functions on $I$ and to $h_I$, since $\int_0^1u(y)\,dy=0$ and $u(1-y)=u(y)$. If $R\ne R'$ in $\cD_m$
intersect, then $I_{\ell}\subsetneq I'_{\ell}$ for some $\ell$; so $I_{\ell}$ lies in a half of $I'_{\ell}$, where
$\Lambda_{I'_{\ell}}$ is affine and $\psi_{I'_{\ell}}$ is constant or $\pm h_{I_{\ell}}$. Hence
$\int\Lambda_{R'}\psi_R\,d\lambda=8^{-d}|R|^2\delta_{RR'}$ and
$\int\psi_R\psi_{R'}\,d\lambda=|R|\delta_{RR'}$.
So $\|g_\eta\|_2^2=(m+1)^{-(d-1)}K_m2^{-m}\le1$. Moreover, $\int g_\eta\,d\lambda=0$, so
$\rho_\eta\ge1-\int(g_\eta)_+\,d\lambda\ge\frac12$, $\int\rho_\eta\,d\lambda=1$ and $\|\rho_\eta\|_2\le2$.

Let $\eta$ be antisymmetric. Since $u(1-y)=u(y)$, we have $\psi_I(1-y)=\psi_{S(I)}(y)$, so
$\psi_R\circ S=\psi_{S(R)}$ and $g_\eta\circ S=\sum_R\eta_R\psi_{S(R)}=-g_\eta$. Hence
$(g_\eta)_+\circ S=(g_\eta)_-$, and as $\int(g_\eta)_+\,d\lambda=\int(g_\eta)_-\,d\lambda$, we get
$\rho_\eta-\rho_\eta\circ S=g_\eta$. Let $\eta'$ be antisymmetric and differ from $\eta$ in at
least $K_m/4$ entries, and let $x_R=\int\Lambda_R(\rho_\eta-\rho_{\eta'})\,d\lambda$. As
$\Lambda_{S(R)}=\Lambda_R\circ S$ and $S$ preserves $\lambda$,
$x_R-x_{S(R)}=\int\Lambda_R(g_\eta-g_{\eta'})\,d\lambda=8^{-d}|R|^2(m+1)^{-\frac{d-1}2}(\eta_R-\eta'_R)$.
Since $(x_R-x_{S(R)})^2\le2x_R^2+2x_{S(R)}^2$, Lemma~\ref{lem:haar} and Bessel's inequality give
\begin{align*}
\|T\rho_\eta-T\rho_{\eta'}\|_2^2\ &\ge\ 2^m\sum_{R\in\cD_m}x_R^2\ \ge\
2^{m-2}\sum_{R\in\cD_m}(x_R-x_{S(R)})^2\\
&\ge\ 2^{m-2}\,\frac{K_m}4\,\frac{4\cdot2^{-4m}}{64^d(m+1)^{d-1}}
\ \ge\ c\,2^{-2m}=:s^2 .
\end{align*}
Lemmas~\ref{lem:few} and~\ref{lem:count} apply with $p=2$, $\xi=1$, $\nu=0$, $\tau=s/6$, the set
$\cG$ of probability densities with $\|\rho\|_2\le2$, $M=16dn^4$, and the $\rho_\eta$ for
antisymmetric $\eta$. Indeed,
$\int_{U_t}\rho\,d\lambda\le\|\rho\|_2(d/M)^{1/2}\le\frac12n^{-2}$ for $t\in D$. For
large $n$, $\log_2\#\cC\le2^{d+4}n\log_2n\le K_m/24$, since $n\le b_d'2^mm^{d-2}$ and
$\log_2n\le2m$. By the choice of $m$, $s\ge c(\log n)^{d-2}/n$, so $n^{-2}\le\tau\le2$ for large
$n$, and $K_m\ge cn\log n$ gives the bound for the proportion.
\end{proof}

\subsection{\texorpdfstring{Upper bound for densities in $L_2$}{Upper bound for densities in L2}}\label{sec:L2dens}

For probability densities with $\|\rho\|_2\le2$, the exponent of the worst case of the
$L_2$-discrepancy is at least $d-2$ by Proposition~\ref{prop:dens}. It is at most $d-1$, even
with nodes that do not depend on the density. This follows from a known bound for
$L_2$-approximation with function values, of Dolbeault, Krieg and Ullrich \cite{DKU} and
Chkifa, Dolbeault, Krieg and Ullrich \cite{CDKU}.

\begin{theorem}[Upper bound for densities in $L_2$]\label{thm:upper}
Let $d\ge2$. There is $C_d<\infty$ such that for every $n\ge2$ there are nodes
$x_1,\dots,x_n\in D$ with the following property. For every probability density
$\rho\in L_2(D)$, there are weights $w_1,\dots,w_n\in\RR$ such that the rule $Q_n$ with these
nodes and weights satisfies
\[
D_2(Q_n;\rho)\ \le\ C_d\,\|\rho\|_2\,\frac{(\log n)^{d-1}}n .
\]
\end{theorem}
\begin{proof}
Let $F_2$ be the set of functions $f(x)=\int_{[x,1)}g\,d\lambda$ with $\|g\|_2\le1$. By
Fubini's theorem, $\int f\rho\,d\lambda-Q_n(f)=\int(T\rho-TQ_n)g\,d\lambda$ for every rule $Q_n$,
so $D_2(Q_n;\rho)=\sup_{f\in F_2}|\int f\rho\,d\lambda-Q_n(f)|$. The set $F_2$ is the unit ball of the reproducing kernel Hilbert space $H$ of these functions, with the norm
$\|f\|_H=\|g\|_2$ and the kernel
$K(x,y)=\int\1_{[x,1)}\1_{[y,1)}\,d\lambda=\prod_{\ell}(1-\max\{x_{\ell},y_{\ell}\})$. The singular values of
the embedding $H\hookrightarrow L_2$ are those of the operator $T^*$ of
Section~\ref{sec:L2thm}, as an operator on $L_2$. In one variable, $T^*$ maps the function
$y\mapsto\sqrt2\sin((j-\frac12)\pi y)$ to $y\mapsto((j-\frac12)\pi)^{-1}\sqrt2\cos((j-\frac12)\pi y)$,
$j\in\NN$, and both systems are orthonormal
bases of $L_2[0,1]$. So the singular values are the non-increasing rearrangement
$\sigma_0\ge\sigma_1\ge\dots$ of the numbers $\prod_{\ell}((j_{\ell}-\frac12)\pi)^{-1}\le\prod_{\ell}j_{\ell}^{-1}$,
$j\in\NN^d$. Since $\#\{j\in\NN^d\colon\prod_{\ell}j_{\ell}\le N\}\le N(1+\log N)^{d-1}$ for $N\ge1$ (by induction on
$d$), we get
$\sigma_k\lesssim k^{-1}(\log(k+1))^{d-1}$ and
$\sum_{k>n/2}\sigma_k^2\lesssim n^{-1}(\log n)^{2d-2}$ for $k\ge1$ and $n\ge2$.
The kernel is bounded and the embedding is injective, so \cite[Corollary~12]{CDKU}, with
$m=\lfloor n/2\rfloor$ there, gives nodes $x_i$ and functions $\phi_i\in L_2(D)$ with
\[
\sup_{f\in F_2}\Bigl\|f-\sum_{i=1}^nf(x_i)\phi_i\Bigr\|_2\ \lesssim\
\sigma_{\lfloor n/2\rfloor}+\Bigl(\frac1n\sum_{k>n/2}\sigma_k^2\Bigr)^{1/2}\ \lesssim\
\frac{(\log n)^{d-1}}n .
\]
See also
\cite[Theorem~3.24 and (3.43)]{KU}; the bound goes back to \cite{DKU}, there with a constant
factor more nodes. Finally, let $w_i=\int_D\phi_i\rho\,d\lambda$. By the
Cauchy--Schwarz inequality,
\[
D_2(Q_n;\rho)=\sup_{f\in F_2}\Bigl|\int_D\Bigl(f-\sum_if(x_i)\phi_i\Bigr)\rho\,d\lambda\Bigr|\le
C_d\|\rho\|_2\frac{(\log n)^{d-1}}n.\qedhere
\]
\end{proof}

\section{Remarks}\label{sec:rem}

\noindent\textbf{Numerical integration.} For $1\le p\le\infty$, let $F_p$ be the set of functions
$f(x)=\nu([x,1))$, $x\in D$, where $|\nu|(D)\le1$ if $p=\infty$, and $\nu$ has a density $g$
with $\|g\|_{p'}\le1$ if $p<\infty$, $\frac1p+\frac1{p'}=1$. Since
$\int f\,d\mu-Q_n(f)=\int(T\mu-TQ_n)\,d\nu$ by Fubini's theorem, the worst-case error of $Q_n$ for
$\int f\,d\mu$ on $F_p$ is $D_p(Q_n;\mu)$ for $p<\infty$ and $D^*(Q_n;\mu)$ for $p=\infty$; this is a
form of the Koksma--Hlawka inequality for general measures \cite{Goetz,AD15}. The class $F_\infty$
consists of functions of bounded variation in the sense of Hardy and Krause \cite{AD15}, and
$F_p$, $p<\infty$, is the unit ball of the Sobolev space of mixed smoothness one with
integrability $p'$ and boundary values zero at the faces $x_{\ell}=1$. So Theorem~\ref{thm:main} is a
lower bound for the worst-case error on $F_\infty$, for every algorithm that uses $n$ function
values, since linear algorithms are optimal for linear functionals on convex and balanced
classes \cite[Chapter~4]{NW}.

\medskip\noindent\textbf{Explicit measures and measures with a density for Theorem~\ref{thm:main}.} Most antisymmetric
sign patterns give hard measures, but we do not know an explicit one. The theorem also holds for measures with a
density: if each atom $x$ of $\mu_\eta$ is spread
uniformly over the cube $x+[0,2^{-3m})^d$, each tent mass changes by at most $d2^{-3m}$, since
$|\Lambda_R(x)-\Lambda_R(y)|\le\sum_{\ell}|x_{\ell}-y_{\ell}|$. This is small compared with
$c_\eta\ge2\theta_m2^{-2m}$ for large $m$, so Lemma~\ref{lem:sep} holds for the smoothed
measures with a smaller constant, and the proof of Theorem~\ref{thm:main}, with
Lemma~\ref{lem:few} applied by the remark after it with $B=n^{3d}$ and $M\ge2dn^{3d+2}$,
and with other constants, gives hard measures with a density bounded
by $2^{3dm}\le n^{3d}$, as $2^m\le n$ by the choice of $m$.

\medskip\noindent\textbf{Symmetric and antisymmetric part.} Let $\mu$ be a probability
measure on $(0,1)^d$, let $\mu'=\mu\circ S^{-1}$ be its mirror image, and let
$\mu_s=\frac12(\mu+\mu')$ and $\mu_a=\frac12(\mu-\mu')$. The proofs of Theorem~\ref{thm:main} and
Corollary~\ref{cor:one} use only $\mu_a$: the tent masses of $\mu$ in $R$ and $S(R)$ differ by
$2\int\Lambda_R\,d\mu_a$. The symmetric part is only needed to make the measure nonnegative,
which holds if and only if $\mu_s\ge|\mu_a|$. For the measures $\mu_\eta$ of Lemma~\ref{lem:loc},
we have $\mu_s=|\mu_a|$, since no atom of $\mu_\eta$ is the mirror image of an atom. So the
method gives nothing for measures that are invariant under $S$, such as the Lebesgue measure.

\medskip\noindent\textbf{$\eps$-approximations.} In computational geometry, a
small set that approximates a finite point set with respect to all ranges of a family is called
an $\eps$-approximation; usually it is a subset of the given set \cite{MV}. For anchored boxes, there are such subsets of size
$O(\eps^{-1}(\log\frac1\eps)^{d-1/2})$ \cite{ABN,MV,Nik}. By Theorem~\ref{thm:main}, for
$d\ge3$, some weighted point sets need at least
$c\,\eps^{-1}(\log\frac1\eps)^{d-1}(\log\log\frac1\eps)^{-\chi_d}$ points, even if these are
arbitrary and carry real weights.

\bigskip\noindent\textbf{Acknowledgements.} This research was funded in whole or in part by the Austrian Science Fund (FWF)
10.55776/PAT3265626. It was also supported by the Austrian Research Promotion Agency (FFG)
through the project FO999921407 (HDcode) funded by the European Union via NextGenerationEU. AI
tools were used for all aspects of this work, including literature search, proofs and writing;
all results and references were carefully checked and revised by the author.

\end{document}